\documentclass{article}
\usepackage[leqno,intlimits]{amsmath}
\usepackage{amssymb}
\usepackage{amsthm}
\usepackage{hyperref}

\newtheorem{tm}{Theorem}[section]
\newtheorem{lm}[tm]{Lemma}
\newtheorem{pr}[tm]{Proposition}
\newtheorem{cor}[tm]{Corollary}

\newcommand*{\wh}{\widehat}
\newcommand*{\un}[1]{\underline{#1}}
\newcommand*{\vp}{\varphi}
\newcommand*{\Cov}{{\text{\bf Cov}}}
\newcommand*{\Zb}{\mathbb Z}
\newcommand*{\Rb}{\mathbb R}
\newcommand*{\te}{\theta}
\newcommand*{\ba}{\begin{aligned}}
\newcommand*{\ea}{\end{aligned}}
\newcommand*{\be}{\begin{equation}}
\newcommand*{\ee}{\end{equation}}
\newcommand*{\e}[1]{\text{\rm e}^{#1}}
\newcommand*{\vr}{\varrho}
\newcommand*{\Ev}{{\bf E}}
\newcommand*{\Pv}{{\bf P}}
\newcommand*{\Vv}{{\text{\bf Var}}}
\newcommand*{\di}{\,\text{\rm d}}
\newcommand*{\wt}{\widetilde}
\newcommand*{\hop}{\bigskip\noindent}
\newcommand*{\xmin}{x^{\text{min}}}
\newcommand*{\xmax}{x^{\text{max}}}
\newcommand*{\musupp}{\mathbb S}  %%support of discrete \mu 
\newcommand*{\nn}{\nonumber} 
\newcommand*{\ve}{\varepsilon} 

\numberwithin{equation}{section}

\begin{document}
\title{A convexity property of expectations under exponential weights}
\author{
\begin{tabular}{@{}c@{\ \,}}
M.\ Bal\'azs\thanks{{\bf Budapest University of Technology and Economics,} Institute of Mathematics, 1 Egry J\'ozsef u., H \'ep.\ V.7., Budapest, Hungary\newline
M. Bal\'azs was partially supported by the Hungarian Scientific Research Fund (OTKA) grants K60708, TS49835, F67729, and the Bolyai Scholarship of the Hungarian Academy of Sciences.
}\\
{\small\tt balazs@math.bme.hu}
\end{tabular}
\begin{tabular}{@{\ \,}c@{}}
T.\ Sepp\"al\"ainen\thanks{{\bf University of Wisconsin-Madison,} Mathematics Department, Van Vleck Hall, 480 Lincoln Dr, Madison WI 53706-1388, USA.
\newline
T.\ Sepp\"al\"ainen was partially supported by National Science Foundation
grant DMS-0402231.} \\
{\small\tt seppalai@math.wisc.edu}
\end{tabular}
}

\maketitle
\begin{abstract}

Please note: after completion of this manuscript we learned that our main results, Theorem \ref{tm:1} and \ref{tm:2}, can be obtained as a special case of Proposition 3.2 on Page 23 of Karlin's book \cite{totpos}.

\hop
Take a random variable \(X\) with distribution \(\mu\) with 
some finite exponential moments, and weight it by an 
exponential factor \(\e{\te X}\) to get the distribution \(\mu^\te\) 
for the admissible  \(\te\)-values.  Define also the 
so-weighted expectation \(\vr(\te):\,=\Ev^\te X\)
with inverse function  \(\te(\vr)\). 
This note proves that for a convex function \(\Phi\), 
\(\Ev^{\te(\vr)}\Phi(X)\) is a convex function of \(\vr\),
 wherever it exists and is finite. Along the way we develop correlation inequalities for convex
functions. Motivation for this result comes from equilibrium investigations of some stochastic interacting systems with
stationary product distributions. In particular, convexity of the
 hydrodynamic flux function follows in some cases.
\end{abstract}


\hop

\section{Introduction}

Please note: after completion of this manuscript we learned that our main results, Theorem \ref{tm:1} and \ref{tm:2}, can be obtained as a special case of Proposition 3.2 on Page 23 of Karlin's book \cite{totpos}.

\hop
Take a non-degenerate random variable \(X\) such that 
\(\Ev\e{\te X}<\infty\) for \(\te\in I\), for some open 
interval $I$. 
For these \(\te\) define the exponentially weighted 
distribution 
\(\Ev^\te(Y):\,=\Ev^\te(Y\e{\te X})/\Ev\e{\te X}\). The exponentially weighted expectation of \(X\) is \(\vr(\te):\,=\Ev^\te X\). 
This function is strictly increasing due to the nondegeneracy assumption.
We denote its inverse by \(\te(\vr)\).

Let \(\Phi\) be a convex function for which
 \(\Ev^{\te(\vr)}\Phi(X)\) exists in an open interval
of \(\vr\)-values.  The first
result of this note is the convexity of the function
\[
\vr\mapsto\Ev^{\te(\vr)}\Phi(X).
\]
Motivation for this  result comes 
from  a class of asymmetric stochastic 
interacting systems that includes the zero range
process (ZRP) and the bricklayer process (BLP). 
We explain these informally 
before turning to precise statements.  

The main application is related to the study of fluctuations of the 
current of particles \(J^{(V)}(t)\)  as seen by an observer moving at a fixed
speed \(V\). Equivalently, these are fluctuations of the height
in the deposition formulation of the process. 
A key fact that underlies some of this work is that the 
variance of the current of a stationary process
 is linked to the deviations of a second
class particle:  
\be
\Vv(J^{(V)}(t))=C(\vr) \widehat\Ev\lvert Q(t)-[Vt]\rvert.
\label{eq:JQ}\ee
The variance on the left is taken in the stationary process
at a fixed density $\vr$. The expectation on the right is
taken with an altered initial product distribution:  the 
density $\vr$ invariant factor is put at each site other than the 
origin, while at the origin there is a different measure 
we denote by \(\nu^{\te(\vr)}\), defined in \eqref{eq:defnu} below.  
$Q(t)$ is the position of a second class particle in the system.
$[Vt]$ is not the usual integer part but rather the integer
between $Vt$ and the origin that is closest to $Vt$. 

Identity \eqref{eq:JQ} has been known for the totally asymmetric simple 
exclusion process since the pioneering work of Ferrari and Fontes
\cite{se}. It was recently extended to the broader class
of processes in \cite{varj2nd}.  For a complete discussion we
 refer the reader  to \cite{varj2nd}, where the measures
  \(\nu^{\te(\vr)}\)  are defined in  equation (2.6) and denoted by  
 \(\wh\mu_\te\).  

The  presently relevant point is that 
for coupling purposes it is important that the measures
  \(\nu^{\te(\vr)}\) are stochastically monotone in the parameter
$\vr$.   In the case of the asymmetric simple exclusion
process (ASEP) this is immediately obvious. This fact was utilized
in a recent coupling-intensive
 proof \cite{se2/3} that established the order $t^{2/3}$
for the variance in \eqref{eq:JQ} for ASEP when the observer
travels at the characteristic speed.  
The first step in a program to extend the variance bounds to
ZRP and BLP is to develop the coupling framework. 
 However, the stochastic monotonicity of the measures
  \(\nu^{\te(\vr)}\) is not at all obvious for ZRP and BLP.
 This we derive from the main convexity
result in Section \ref{sc:mon}. 

The conserved quantity in these asymmetric processes (typically viewed
as particle counts, but also discrete gradients of the 
interface height) satisfies a hydrodynamic scaling limit 
where the limiting evolution is the entropy
solution of a scalar conservation law
of the form 
\[
\partial_t\vr+c\partial_x\mathcal H(\vr)=0.
\] 
(We refer the reader to \cite{cl} for general theory.)
The key quantity is the flux $\mathcal H(\vr)$ which is computed
as the expected  jump rate in the stationary process at
particle density $\vr$.  The constant $c$ is the mean increment
of a particle that has decided to jump. 

A second application of the main convexity result
 is to give an alternative
  proof of the convexity of the hydrodynamic 
flux function $\mathcal H$
for zero range and bricklayer processes when the jump rate is convex. (See \cite{fluct} for the 
original proof.) Concavity of the hydrodynamic flux 
also follows for concave jump rates in the zero range process.
Strict convexity or concavity are also discussed.

The characteristic speed referred to above is $V^\vr=c\mathcal H'(\vr)$.
So the issue can also be framed as the monotonicity 
of this quantity in the particle
density. 

The rest of this note is organized as follows. 
We rewrite the convexity problem in terms of correlation inequalities 
of functions of \(X\).  These we handle via separation of positive and 
negative parts and further correlation inequalities. This is done in Section \ref{sc:comp}.
In Section \ref{sc:ips} we derive 
the consequences  for  stochastic interacting systems.

\section{Derivatives and correlations}\label{sc:comp}

As in the introduction, let $X$ be a nondegenerate real-valued
 random variable
and $\Phi$ a convex function on some 
interval that contains the range of $X$.  The standing assumption
throughout this section  is that
for some open interval $I\subseteq\Rb$, 
\be
\Ev(\e{\te X})<\infty
\quad\text{and}\quad
\Ev(X^2 \lvert\Phi(X)\rvert\e{\te X})<\infty
\label{eq:ass}\ee
for all $\te\in I$.  Define the exponentially weighted distribution
as \(\Ev^\te(Y)=\Ev^\te(Y\e{\te X})/\Ev\e{\te X}\).
The function
 \(\vr(\te)=\Ev^\te X\) is strictly increasing (justification below
in Corollary \ref{cr:monv}). It has an  
 inverse function \(\te(\vr)\) defined in some nontrivial open interval $J$. 
The expectation  $\Ev^{\te(\vr)}\Phi(X)$ is well defined for $\vr\in J$.
 
\begin{tm}
The function $\vr\mapsto \Ev^{\te(\vr)}\Phi(X)$ is convex on $J$.  
\label{tm:1}\end{tm}

Our main interest lies in discrete distributions so we state a
further condition for strict convexity for that case.
Suppose the distribution $\mu$ of $X$ is supported on a discrete
subset $\musupp$  of $\Rb$.  So $\musupp$ is either
 finite or countably infinite but locally finite.   
Let $\Phi$ be a function defined
on $\musupp$.  Extend $\Phi$ to a function
on the smallest closed interval that contains $\musupp$ 
 by connecting adjacent points on the graph of $\Phi$
with line segments. Assume the function $\Phi$ thus defined is
convex.  Say $\Phi$ is {\sl strictly convex} at a point 
$z\in\musupp$ if the slope of the extended $\Phi$ jumps at $z$. Such a 
point $z$ cannot be the maximum or minimum of $\musupp$ because
we have not defined $\Phi$ outside the smallest interval 
that contains $\musupp$.  If no such point exists then $\Phi$ 
is linear. 

\begin{tm}  Consider the discrete case described in the paragraph
above.  
The function $\vr\mapsto \Ev^{\te(\vr)}\Phi(X)$ is
strictly convex throughout its interval of definition if and only if   
$\musupp$ contains  at least three points and $\Phi$ is 
strictly convex at some point of $\musupp$.  In the complementary 
case the function  $\vr\mapsto \Ev^{\te(\vr)}\Phi(X)$ is linear.  
\label{tm:2}\end{tm}

The remainder of this section covers the proofs. Throughout
we only consider values $\te\in I$ for which assumption \eqref{eq:ass} 
guarantees that the 
derivatives and other operations we perform are justified. 
In particular, since $\lvert X\rvert^k\le k!\ve^{-k}(\e{\ve X}+\e{-\ve X})$,
$X$ has all moments under $\Ev^\te$ for each $\te\in I$ 
because $\te\pm\ve\in I$ for small enough $\ve>0$.  

We start with a preliminary lemma, repeated from 
 Lemma A2 of \cite{exists}.
\begin{lm}
For any function \(\vp\), we have
\[
\frac{\di}{\di\te}\,\Ev^\te\vp(X)=\Cov^\te(\vp(X),\,X)
\]
 provided the expectations exist in a neighborhood of \(\te\).
\end{lm}
\begin{proof}  
\[
\ba
\frac{\di}{\di\te}\,\Ev^\te\vp(X)&=\frac{\di}{\di\te}\,\frac{\Ev(\vp(X)\cdot\e{\te X})}{\Ev\e{\te X}}\\
&=\frac{\Ev(\vp(X)\cdot X\cdot\e{\te X})}{\Ev\e{\te X}}-\Ev(\vp(X)\cdot\e{\te X})\cdot\frac{\Ev(X\cdot\e{\te X})}{[\Ev\e{\te X}]^2}\\
&=\Cov^\te(\vp(X),\,X).
\ea
\]
\end{proof}
Recall that we exclude the degenerate
 case where \(\mu\) is supported  on a single point.
\begin{cor}\label{cr:monv}
\begin{align}
\frac{\di\vr(\te)}{\di\te}&=\frac{\di}{\di\te}\,\Ev^\te X=
\Cov^\te(X,\,X)=\Vv^\te X>0\text{,\qquad and}\nn\\
\frac{\di\Ev^{\te(\vr)}\vp(X)}{\di\vr}
&=\frac{\di\Ev^{\te}\vp(X)}{\di\te}\cdot\frac{\di\te(\vr)}{\di\vr}=\frac{\Cov^{\te(\vr)}(\vp(X),\,X)}{\Vv^{\te(\vr)}X}.
\label{eq:ddvr}\end{align} 
\end{cor}

\noindent
Now we proceed by rewriting the second derivative of \(\Ev^{\te(\vr)}\Phi(X)\) in terms of covariances. We omit the notation \((X)\) from \(\Phi(X)\).
\begin{lm}\label{lm:eq}
The following are equivalent:
\renewcommand{\theenumi}{\alph{enumi})}
\renewcommand{\labelenumi}{\alph{enumi})}
\begin{enumerate}
\item For any convex function \(\Phi\), \(\Ev^{\te(\vr)}\Phi\) is a convex function of \(\vr\).
\item For any convex function \(\Phi\),
\be
\Cov^\te(\wt\Phi\cdot X,\,X)\cdot\Cov^\te(X,\,X)\ge\Cov^\te(\Phi,\,X)\cdot\Cov^\te(\wt X\cdot X,\,X),\label{eq:b}
\ee
where \(\wt\cdot\) stands for centering w.r.t.\ \(\Ev^\te\).\label{eq1b}
\end{enumerate}
\end{lm}
\begin{proof}
We write, as in the corollary,
\[
\frac{\di}{\di\vr}\Ev^{\te(\vr)}\Phi=\frac{\Cov^{\te(\vr)}(\Phi,\,X)}{\Vv^{\te(\vr)}X}.
\]
We need to see if this is nondecreasing in \(\vr\)
or, equivalently,  nondecreasing in \(\te\). That happens if and only if
\[
\ba
0&\le\Vv^\te X\cdot\frac{\di}{\di\te}\Cov^\te(\Phi,\,X)-\Cov^\te(\Phi,\,X)\cdot\frac{\di}{\di\te}\Vv^\te X\\
&=\Vv^\te X\cdot\bigl[\Cov^\te(\Phi X,\,X)-\Cov^\te(\Phi,\,X)\cdot\Ev^\te(X)-\Ev^\te(\Phi)\cdot\Cov^\te(X,\,X)\bigr]\\
&\quad-\Cov^\te(\Phi,\,X)\cdot[\Cov^\te(X^2,\,X)-2\Ev^\te X\cdot\Cov^\te(X,\,X)]\\
&=\Cov^\te(\wt\Phi X,\,X)\cdot\Cov^\te(X,\,X)-
\Cov^\te(\Phi,\,X)\cdot\Cov^\te(\wt XX,\,X).
\qedhere\ea
\]
\end{proof}

\noindent
Next we concentrate on proving that part \ref{eq1b} of the
last lemma holds for
{\sl any distribution}. Therefore we omit the superscript \(\te\).
\begin{lm}
Part \ref{eq1b} of Lemma \ref{lm:eq} is further equivalent to each of these
two statements:
\begin{itemize}
\item[c)] For any convex function \(\Phi\) that is uncorrelated with \(X\), 
\(\Cov(\Phi,\,X^2)\ge0\).
\item[d)] For any convex function \(\Phi\),
\be
\Cov(\Phi,\,X^2)\cdot\Cov(X,\,X)\ge\Cov(\Phi,\,X)\cdot\Cov(X^2,\,X).\label{eq:d}
\ee
\end{itemize}
\end{lm}
\begin{proof}
Given a convex function \(\Phi\), let 
\(\wh\Phi(X):\,=\Phi(X)-C\cdot X\) with $C$ chosen so that
\(\wh\Phi\) is uncorrelated with \(X\).
 \(\wh\Phi\) is also  convex, and we note 
 that \eqref{eq:b} holds for \(\Phi\) 
if and only if it holds for \(\wh\Phi\). Hence \ref{eq1b} is equivalent to the statement obtained by restricting \ref{eq1b} to convex functions that are uncorrelated with \(X\). For such functions this statement becomes
\[
\ba
0&\le\Cov(\wh\Phi X,\,X)-\Ev\wh\Phi\cdot\Cov(X,\,X)\\
&=\Cov(\wh\Phi,\,X^2)+\Ev\wh\Phi\cdot\Ev(X^2)-\Ev(\wh\Phi X)\cdot\Ev X-\Ev\wh\Phi\cdot\Ev(X^2)+\Ev\wh\Phi\cdot\Ev X\cdot\Ev X\\
&=\Cov(\wh\Phi,\,X^2).
\ea
\]
Thus b) is equivalent to c). 

Condition c) is a weakening of d), and we see that
 c) implies d) by determining the constant in the transformation 
that led to \(\wh\Phi\):
\[
\Cov(\wh\Phi,\,X)=\Cov(\Phi,\,X)-C\cdot\Cov(X,\,X)=0,
\]
therefore
\[
\wh\Phi=\Phi-\frac{\Cov(\Phi,\,X)}{\Cov(X,\,X)}\cdot X.
\]
Substituting  this into \(\Cov(\wh\Phi,\,X^2)\ge0\) of c) leads to d).
\end{proof}

Next we show that for part d) of the above lemma it suffices to 
consider the special case  \(\Phi(X)=|X|\).

\begin{lm}\label{lm:abs}
Part d) of the above lemma is implied by this statement: 
\begin{itemize}
\item[e)] For any distribution (with finite third absolute moments) we have
\be
\Cov(|X|,\,X^2)\cdot\Cov(X,\,X)\ge\Cov(|X|,\,X)\cdot\Cov(X^2,\,X).\label{eq:e}
\ee
\end{itemize}
\end{lm}
\begin{proof}
Consider functions of the form
\be
\phi(x)=c+ a\cdot[x-x_0]^+-b\cdot[x-x_0]^-\label{eq:form}
\ee
for some \(a> b\) and \(x_0,\,c\in\Rb\). Notations \(^+\) and \(^-\) stand
 for positive and negative parts, respectively. These functions are convex.
The first claim is  that if \eqref{eq:d} holds for functions of this 
special form, then it holds for any convex \(\Phi\).

This follows because $\Phi$ can be approximated from below in a 
pointwise fashion  by a sequence of 
functions of this type: 
\[
g(x) = c -a_0[x-y_1]^-+a_1[x-y_1]^+ + \sum_{k=2}^m (a_k-a_{k-1}) [x-y_k]^+
\]
with $a_0<a_1<\dotsm<a_m$ and $y_1<\dotsm<y_m$. 
The function $g$ above is a sum of convex functions of type \eqref{eq:form}.
To see the approximation, take points $z_0<z_1<\dotsm<z_m$ and 
let $a_i$ be the slope of a tangent to $\Phi$ at the point
 $(z_i, \Phi(z_i))$.  
Pick the $z_i$'s so that the $a_i$'s are strictly increasing. 
(This entails no loss of generality because a linear approximation
to $\Phi$ is exact throughout any interval with constant slope.)  
Let $g_i$ ($0\le i\le m$) be the linear function of slope $a_i$ that passes 
 through the point $(z_i, \Phi(z_i))$.  Let
 $y_i$  ($1\le i\le m$) be the $x$-coordinate
of the point where the graphs of $g_{i-1}$ and $g_i$ intersect and set
\[
c=\Phi(z_0)+ a_0(y_1-z_0)=\Phi(z_1)+ a_1(y_1-z_1).
\]
 Then it can be checked that $g$ from above is the 
pointwise maximum of the $g_i$'s, or equivalently, that 
$g=g_i$ on $(y_{i},y_{i+1})$ with $y_0=-\infty$
and $y_{m+1}=\infty$.  
By choosing the $z_i$'s
carefully one can create a sequence of convex functions $g^{(m)}$ such that 
 $g^{(m)}\nearrow\Phi$ pointwise.
By \eqref{eq:ass}
monotone convergence applies to show
 \(\Cov(g^{(m)},\,X^b)\to \Cov(\Phi,\,X^b)\) for $b=1,2$.

Thus we can derive  \eqref{eq:d} for $\Phi$ by checking it
for each $g^{(m)}$.  
Since \eqref{eq:d} is linear in \(\Phi\), it is then enough to know that it holds 
for each term of the type \eqref{eq:form}. This we now check. 

\hop

With suitably chosen constants \(A>0\), \(B\) and \(C\), the transformation
\be
\phi(x)\mapsto A\phi(x)+Bx+C
\label{eq:trans}\ee
turns \(\phi\) of \eqref{eq:form} into the function \(|x-x_0|\).
(Note that $a>b$ is needed for this.)
 The left and right-hand sides of  \eqref{eq:d} are, up to the multiplying factor \(A\),
 invariant under these transformations.  
Hence \eqref{eq:d}   holds for \(\phi\) if and only if it holds for \(|x-x_0|\):
\[
\Cov(|X-x_0|,\,X^2)\cdot\Cov(X,\,X)\ge\Cov(|X-x_0|,\,X)\cdot\Cov(X^2,\,X).
\]
Introduce now \(Y=X-x_0\), and write this inequality in the form
\[
\Cov(|Y|,\,(Y+x_0)^2)\cdot\Cov(Y,\,Y)\ge\Cov(|Y|,\,Y)\cdot\Cov((Y+x_0)^2,\,Y).
\]
Subtracting \(2x_0\cdot\Cov(|Y|,\,Y)\cdot\Cov(Y,\,Y)\) from both sides leads
 to e) (for the distribution of \(Y=X-x_0\)).
\end{proof}

\noindent
Some elementary computations will now finish the proof of Theorem \ref{tm:1}.
\begin{lm}\label{lm:pn}
Part e) in Lemma \ref{lm:abs} holds.
\end{lm}
\begin{proof}
For this proof, we introduce the positive and negative part moments:
\[
P_i:\,=\Ev\bigl((X^+)^i\bigr),\qquad N_i:\,=\Ev\bigl((X^-)^i\bigr).
\]
Expanding \eqref{eq:e} gives
\begin{multline*}
\bigl[P_3+N_3-(P_1+N_1)\cdot(P_2+N_2)\bigr]\cdot\bigl[P_2+N_2-(P_1-N_1)^2\bigr]\\
\ge\bigl[P_2-N_2-(P_1+N_1)\cdot(P_1-N_1)\bigr]\cdot\bigl[P_3-N_3-(P_2+N_2)\cdot(P_1-N_1)\bigr].
\end{multline*}
Somewhat tedious factoring shows that this is equivalent to
\begin{align}
0&\le N_1\cdot(P_3P_1-P_2P_2)\label{eq:1st}\\
&\quad+P_1\cdot(N_3N_1-N_2N_2)\label{eq:2nd}\\
&\quad+P_2N_3-P_1P_1N_3-P_2N_2N_1\label{eq:3rd}\\
&\quad+P_3N_2-P_3N_1N_1-P_2P_1N_2.\label{eq:4th}
\end{align}
We proceed by showing that each line above
 is non-negative. Clearly if \(\Pv\{X>0\}\) or \(\Pv\{X\le0\}\) is zero, then \(P_i\)'s or \(N_i\)'s are zero and the statement is trivially true. Assuming the contrary and dividing \eqref{eq:1st} by \([\Pv\{X>0\}]^2\) makes conditional expectations out of the \(P_i\)'s:
\[
\frac{P_3P_1-P_2P_2}{[\Pv\{X>0\}]^2}=\Ev(X^3\,|\,X>0)\cdot\Ev(X\,|\,X>0)-\Ev(X^2\,|\,X>0)\cdot\Ev(X^2\,|\,X>0).
\]
To show that this is non-negative, introduce the expectation
\[
\wh\Ev(\cdot):\,=\frac{\Ev(\cdot\times X^2\,|\,X>0)}{\Ev(X^2\,|\,X>0)},
\]
with which the previous formula becomes a constant multiple of
\[
\wh\Ev X\cdot\wh\Ev\frac1X-1=-\wh\Cov\Bigl(X,\,\frac1X\Bigr).
\]
Notice that the \(\wh{\phantom{a}}\) measure is concentrated on positive values, where \(1/X\) is a decreasing function of \(X\) hence the above covariance is non-positive. A similar argument shows 
that \eqref{eq:2nd} is non-negative.

Separate  \eqref{eq:3rd} into the sum of two terms: 
\[
\bigl[P_2\Pv\{X>0\}-P_1P_1\bigr]\cdot N_3+P_2\cdot\bigl[N_3\Pv\{X\le0\}-N_2N_1\bigr].
\]
Divide the first bracket by \([\Pv\{X>0\}]^2\) to get
\[
\Ev(X^2\,|\,X>0)-[\Ev(X\,|\,X>0)]^2\ge0.
\]
Dividing the second bracket by \([\Pv\{X\le 0\}]^2\) leads to
\be
\begin{split}
&\Ev(|X|^3\,\bigr|\,X\le0)-\Ev(X^2\,\bigr|\,X\le0)\cdot\Ev(|X|\,\bigr|\,X\le0)\\
&\qquad\qquad\qquad
=\Cov(X^2,\,|X|\,\bigr|\,X\le0)\ge0
\end{split}\label{eq:aux6}\ee
since \(X^2\) is an increasing function of \(|X|\) on non-positive numbers. The term \eqref{eq:4th} is treated in a similar manner.
\end{proof}

Tracing the lemmas  backward shows that we have verified part a) of 
Lemma \ref{lm:eq} and thereby proved Theorem \ref{tm:1}.

To prove  Theorem \ref{tm:2}, note first that in the
complementary case $\Phi(X)=aX$ on $\musupp$, and then
\eqref{eq:ddvr} implies that the derivative 
${\di\Ev^{\te(\vr)}\Phi(X)}/{\di\vr}$ 
is constant.  

To prove the main statement of Theorem \ref{tm:2} we retrace
some earlier steps. Let $x_0\in\musupp$ be a point of
strict convexity whose
existence is assumed. Namely, 
\be
\Pv\{X<x_0\}\Pv\{X=x_0\}\Pv\{X>x_0\}>0
\label{eq:aux9}\ee
 and the slopes
$b=\Phi'(x_0-)$ and $a=\Phi'(x_0+)$ satisfy $a>b$.  
Then we can write
\[
\Phi(x)=\Phi(x_0)+a(x-x_0)^+-b(x-x_0)^-+\Psi(x)
\]
for another convex function $\Psi$ that vanishes on an interval
around $x_0$.  Since we already have Theorem \ref{tm:1} 
for $\Psi$,  it suffices to prove strict convexity of
$\vr\mapsto\Ev^{\te(\vr)}\phi(X)$ for 
\[
\phi(x)= \Phi(x_0)+a(x-x_0)^+-b(x-x_0)^-.
\]
After an application of the transformation \eqref{eq:trans} 
the question boils down to showing {\sl strict} inequality 
in \eqref{eq:e} for the new variable $Y=X-x_0$.  For this
it suffices to check that at least one of the quantities
\eqref{eq:1st}--\eqref{eq:4th} is strictly positive. 
From \eqref{eq:aux9} follows that each $P_i$ and $N_i$ is 
strictly positive.  Schwarz inequality shows that 
$P_2P_2<P_3P_1$ if $Y$ has two distinct strictly positive values,
and  
$N_2N_2<N_3N_1$ if $Y$ has two distinct strictly negative values.
If both these requirements fail, then \eqref{eq:aux9}
forces $Y$ to take
one positive value, one negative value, and the value zero
with positive probability.    But then this makes the quantity
in \eqref{eq:aux6} strictly positive for $Y$.  

Thus we conclude that strict inequality holds
in \eqref{eq:e} for $Y=X-x_0$, and strict convexity
of $\vr\mapsto \Ev^{\te(\vr)}\phi(X)$ follows.  We have proved
Theorem \ref{tm:2}. 

\section{Application to stochastic interacting systems}\label{sc:ips}

To keep this note short we give a minimal possible introduction
to the applications of the convexity result and refer the reader
  to \cite{fluct} and \cite{varj2nd} for the complete picture. 
Let \(-\infty\le\xmin\le0\) and \(1\le\xmax\le\infty\) be (possibly infinite valued) integers, and consider the discrete interval 
\(I=(\xmin-1,\,\xmax+1)\cap\Zb\). 
Fix a function \(f\,:\,I\to\Rb^+\). For $I\ni x>0$ we set
\[
f(x)!:\,=\prod_{y=1}^xf(y),
\]
while for $I\ni x<0$ let
\[
f(x)!:\,=\frac{1}{\prod\limits_{y=x+1}^0f(y)},
\]
finally $f(0)!:\,=1$. Then we have 
\[
f(x)!\cdot f(x+1)=f(x+1)!
\]
for all $x\in I$. Let
\[
\bar\te:\,=\left\{\begin{array}{ll}\log\left(\liminf\limits_{x\to\infty}\left(f(x)!\right)^{1/x}\right)\ \ &,\ \text{if}\ \xmax=\infty\\\infty\ \ &,\ \text{else}\end{array}\right.
\]
and
\[
\un\te:\,=\left\{\begin{array}{ll}\log\left(\limsup\limits_{x\to\infty}\left(f(-x)!\right)^{-1/x}\right)\ \ &,\ \text{if}\ \xmin=-\infty\\-\infty\ \ &,\ \text{else}.\end{array}\right.
\]
We require \(f\) to be such that \(\un\te<0<\bar\te\). In this case
\[
\mu(x):\,=\frac{\frac{1}{f(x)!}}{\sum\limits_{y\in I}\frac{1}{f(y)!}}
\]
defines a probability measure on \(I\), and the exponentially weighted version
\[
\mu^\te(x):\,=\frac{\frac{\e{\te x}}{f(x)!}}{\sum\limits_{y\in I}\frac{\e{\te y}}{f(y)!}}
\]
is also well defined for any \(\un\te<\te<\bar\te\). 
This latter is the marginal of a  stationary product 
distribution of many stochastic 
interacting systems, see e.g.\ \cite{varj2nd}. 

\subsection{Convexity of hydrodynamic flux
 for zero range and bricklayer processes}

In particular, the attractive zero range process is an 
example where \(I=[0,\,\infty)\cap\Zb\), \(f(0)=0<f(1)\), and \(f\) is non-decreasing.  The rate for a particle to 
jump from a site with \(x\) particles is \(f(x)\). Its hydrodynamic
(macroscopic)  
flux function \(\mathcal H\,:\,\Rb^+\to\Rb^+\) is given by 
\[
\mathcal H(\vr) =\Ev^{\te(\vr)}f(X)
\]
with the notation of the Introduction. The results of the previous 
section for \(f\)  now read as follows:
\begin{pr}
If the jump rate \(f\) of the zero range process is convex (or concave), then the hydrodynamic flux \(\mathcal H\) is also convex (or concave, respectively). Moreover, in this case \(\mathcal H\) is strictly convex (or concave, respectively) if and only if \(f\) is not linear.
\end{pr}

The bricklayer process has \(I=(-\infty,\,\infty)\cap\Zb\) 
and \(f\) non-decreasing such that \(f(x)\cdot f(1-x)=1\) for all \(x\in\Zb\). Its jump rate for a brick to be laid on a column between negative discrete gradients \(x\) on the left and \(y\) on the right is \(f(x)+f(-y)\), see \cite{varj2nd} for more 
details. The hydrodynamic flux function \(\mathcal H\,:\,\Rb\to\Rb^+\) is now
\[
\mathcal H(\vr) =\Ev^{\te(\vr)}\bigl(f(X)+f(-Y)\bigr)
\]
where \(X\) and \(Y\) are i.i.d.\ variables with distribution 
\(\mu^{\te(\vr)}\). Notice that 
 non-decreasingness and non-negativity of \(f\) on \(\Zb\) 
excludes concave functions with the exception of the constant one function. Our result for this process is
\begin{pr}
If the function \(f\) of the bricklayer process is convex and not constant one, then its hydrodynamic flux \(\mathcal H\) is strictly convex.
\end{pr}
Parts of these two propositions were proved with  coupling 
methods in \cite{fluct}.

\subsection{Monotonicity of a special distribution}
\label{sc:mon}
We come to the primary motivation of the note. 
As explained in the Introduction, the study of current fluctuations
uses couplings of processes whose initial particle number at the
origin obeys the following type of distribution:
\begin{equation}
\nu^{\te(\vr)}(y)=\frac{1}{\Vv^{\te(\vr)} X}\sum_{x=y+1}^{\xmax}[x-\Ev^{\te(\vr)} X]\cdot\mu^{\te(\vr)}(x)\quad(\xmin\le y<\xmax).
\label{eq:defnu}\end{equation}
(See \cite[eqn.~(2.6)]{varj2nd} for the original definition.)
To create couplings with useful monotonicity properties,
one needs these 
 distributions to be  monotone in the 
parameter \(\vr\), in the sense of stochastic domination.  
This we can now derive as a consequence 
of the main result. 
\begin{pr}
The family of measures \(\nu^{\te(\vr)}\) is monotone in \(\vr\).
\end{pr}
\begin{proof}
By Corollary \ref{cr:monv},
\[
\ba
\nu^\te(y)&=\frac{1}{\Vv^\te X}\cdot\Ev\bigl([X-\Ev^\te(X)]\cdot{\bf1}\{X>y\}\bigr)\\
&=\frac{\Cov^\te(X,\,{\bf1}\{X>y\})}{\Vv^\te X}=\frac{\di}{\di\vr}\Pv^{\te(\vr)}\{X>y\}.
\ea
\]
Let us denote the \(\nu^{\te(\vr)}\)-expectation by \(\Ev^{\nu,\,\te(\vr)}\). Monotonicity of the family \(\nu^{\te(\vr)}\) is equivalent to the
 property that, for any bounded non-decreasing function \(\vp\),
\[
0\le\frac{\di}{\di\vr}\Ev^{\nu,\,\te(\vr)}\vp(X).
\]
We compute a different expression for this derivative.  Passing the
derivative through the sum in the third equality below is 
justified because the series involved are dominated by   
certain geometric series, uniformly over $\te$  in small 
open neighborhoods.  This follows from the definitions of
 \(\un\te\) and \(\bar\te\) and the assumption \(\un\te<0<\bar\te\).
\[
\ba
\Ev^{\nu,\,\te(\vr)}\vp(X)&=\sum_{y=\xmin}^{\xmax}\vp(y)\cdot\frac{\di}{\di\vr}\Pv^{\te(\vr)}\{X>y\}\\
&=\sum_{y=\xmin}^{\xmax}\vp(y)\cdot\frac{\di}{\di\vr}[\Pv^{\te(\vr)}\{X>y\}-{\bf1}\{0\ge y\}]\\
&=\frac{\di}{\di\vr}\sum_{y=\xmin}^{\xmax}\vp(y)\cdot[\Pv^{\te(\vr)}\{X>y\}-{\bf1}\{0\ge y\}]\\
&=\frac{\di}{\di\vr}\Ev^{\te(\vr)}\sum_{y=\xmin}^{\xmax}\vp(y)\cdot[{\bf1}\{X>y\}-{\bf1}\{0\ge y\}]\\
&=\frac{\di}{\di\vr}\Ev^{\te(\vr)}\sum_{y=\xmin}^{\xmax}\vp(y)\cdot[{\bf1}\{X>y>0\}-{\bf1}\{0\ge y\ge X\}]\\
&=\frac{\di}{\di\vr}\Ev^{\te(\vr)}\Bigl[\sum_{y=1}^{X-1}\vp(y)-\sum_{y=X}^0\vp(y)\Bigr]=\frac{\di}{\di\vr}\Ev^{\te(\vr)}\Phi(X).
\ea
\]
Above we introduced 
 the function \[\Phi(x)=\sum\limits_{y=1}^{x-1}\vp(y)-
\sum\limits_{y=x}^0\vp(y), \]
 with the convention that empty sums are zero. To conclude the proof, notice that \(\Phi(x+1)-\Phi(x)=\vp(x)\). 
Thus a non-decreasing function \(\vp\) determines a (non-strictly) 
convex function \(\Phi\) with \(\Phi(1)=0\), and vice-versa. 
Hence Section \ref{sc:comp} establishes that
\[
\frac{\di}{\di\vr}\Ev^{\nu,\,\te(\vr)}\vp(X)=\frac{\di^2}{\di\vr^2}\Ev^{\te(\vr)}\Phi(X)\ge0.
\qedhere\]
\end{proof}

\section*{Acknowledgment}

We thank Omer Angel and B\'alint T\'oth for very illuminating discussions on the subject.

\bibliography{refsmarton}
\bibliographystyle{plain}

\end{document}